\documentclass[11pt]{article}

\usepackage{amsfonts,amssymb,wrapfig}
\usepackage[all, knot]{xy}

\newcommand{\R}{\mathbb{R}}

\newtheorem{theorem}{Theorem}

\input{epsf.sty}
\usepackage{graphicx}

\begin{document}

\title{Reidemeister/Roseman-type Moves to Embedded Foams in $4$-dimensional Space}

\author{
J. Scott Carter
\\ University of South Alabama}

\maketitle

\begin{abstract} 

The dual to a tetrahedron consists of a single vertex at which four edges and six faces are incident. Along each edge, three faces converge. 
 A $2$-foam is a compact topological space such that each point has a neighborhood homeomorphic to a neighborhood of that complex. 
 Knotted foams in $4$-dimensional space are to knotted surfaces, as knotted trivalent graphs are to classical knots. The diagram of a knotted foam consists of a generic projection into $3$-space with crossing information indicated via a broken surface. In this paper, a finite set of moves to foams are presented that are analogous to the Reidemeister-type moves for knotted graphs. These moves include the Roseman moves for knotted surfaces. Given a pair of diagrams of isotopic knotted foams there is a finite sequence of moves taken from this set that, when applied to one diagram sequentially, produces the other diagram.

\end{abstract} 

\section{Introduction}

\begin{wrapfigure}{r}{2.7in}\vspace{-.4in}\includegraphics[width=2.5in]{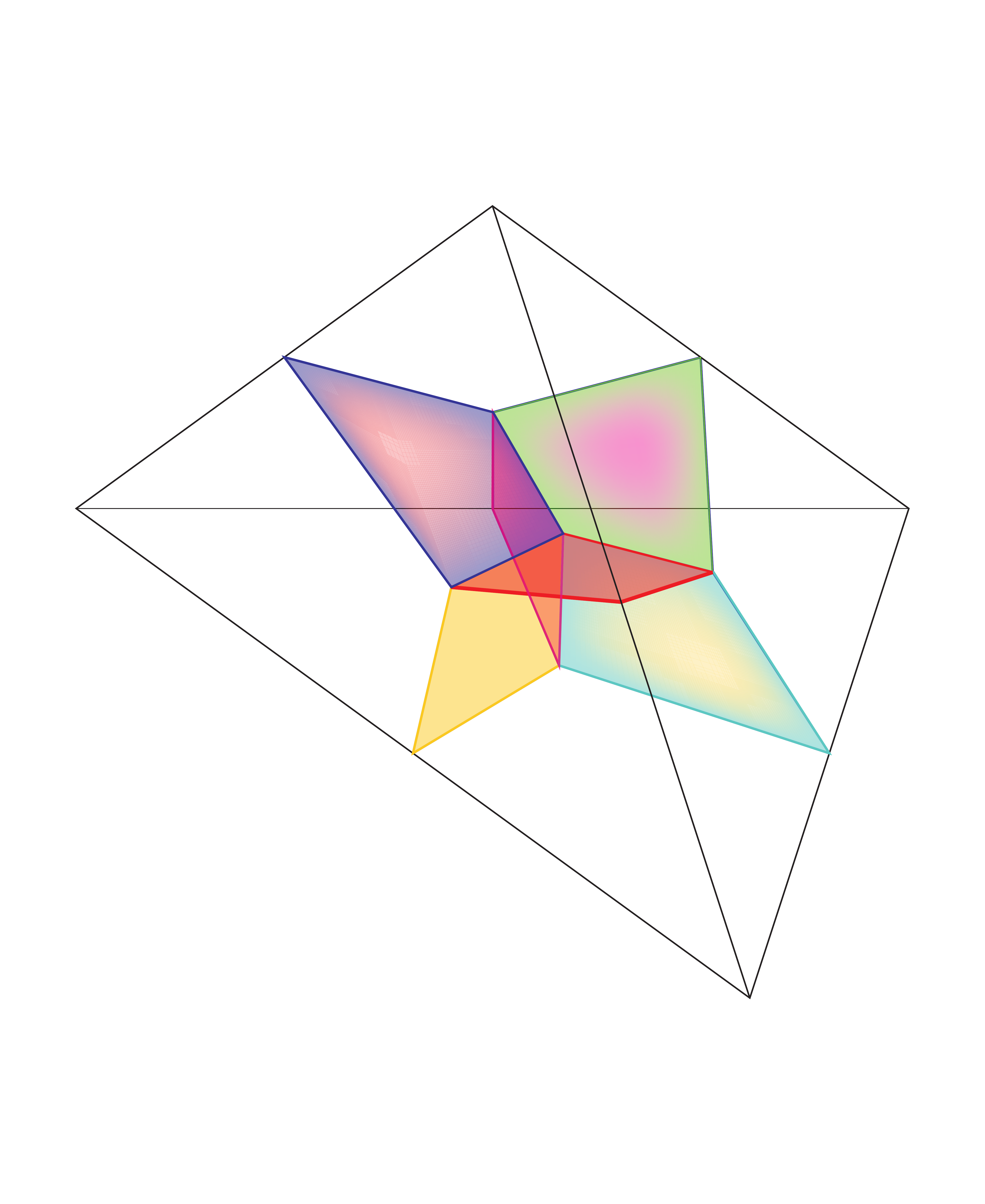}\end{wrapfigure}
Consider the spine of the tetrahedron that is obtained by embedding four copies of the topological space that is homeomorphic to the alpha-numeric character ${\sf Y}$ in each of the triangular faces of the tetrahedron and coning the result to the barycenter. This two dimensional space (illustrated to the right and below), $Y^2$, has a single vertex, four edges, and six $2$-dimensional faces. Three faces are incident to each edge, and a neighborhood of a point in an open edge is homeomorphic to  ${\mbox{\sf Y}}  \times (-1,1).$ A {\it $2$-dimensional foam ($2$-foam)} is a compact topological space, $F$, such that any point has a neighborhood that is homeomorphic to a neighborhood of a point in $Y^2$. Thus a foam is stratified into isolated singular points, $1$-dimensional edges at which three sheets meet, and $2$-dimensional faces. The boundary of a foam is a trivalent graph. 
A {\it closed foam} has empty boundary.
Analogous concepts exist in all dimensions. Just as a trivalent graph can be embedded and knotted in $3$-space, a $2$-foam can be embedded and knotted in $4$-dimensional space.

\begin{wrapfigure}{l}{2.2in}\includegraphics[width=2in]{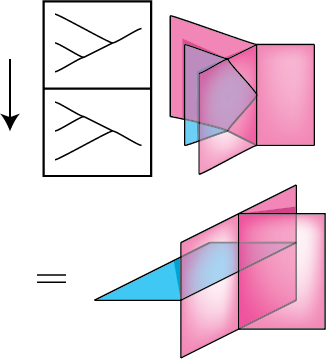}\end{wrapfigure}The space $Y^2$  can be interpreted via a movie of the associativity rule when this is expressed in terms of binary trees.
Foams are important since they include special spines of $3$-dimensional manifolds \cite{Matveev,Piergallini}. They are related to categorifications of the HOMFLYPT (FLYTHOMP) polynomial \cite{Kho,Vaz,MacVaz}.  Knotted closed $2$-foams and their higher dimensional generalizations can be used to represent $3$-cycles in a homology theory of $G$-family of quandles and other more general algebraic structures. 

The purpose of this paper is to present an analogue of the Reidemeister-type moves for knotted foams in $4$-space. Just as knotted trivalent graphs (spacial graphs) contain classical knots and links as a subset, knotted foams include embedded surfaces in $4$-dimensional space. Thus the moves that will be presented will include the Roseman moves \cite{Roseman}. Indeed, the proof that the given set of moves is sufficient to transform two diagrams of isotopic foams follows closely Roseman's original proof of the sufficiency of his set of seven moves. 

To achieve the goals of presenting a set of moves to foams and demonstrating their sufficiency, the local pictures  that are used to describe knottings of foams are given. These local crossings are obtained from the Reidemeister-type moves to graphs. The sufficiency of such moves are obtained by examining the generic critical points and transverse intersections of the self-intersection strata. Indeed, the description of the moves for foams are precisely an analysis of critical behavior and intersections between self-intersections or edges of the foam.

The stage will be set  following these acknowledgement.

\subsection*{Acknowledgements} This paper was studied with the support of the Ministry of Education Science and Technology (MEST) and the Korean Federation of Science and Technology Societies (KOFST). 
I owe much of the approach here to a conversation that I had with Osamu Saeki. This paper is part of an on-going project with Atsushi Ishii and Masahico Saito. In addition, I have had many valuable conversations with Seiichi Kamada, Shin Satoh, and the faculty and students at the TAPU workshops and seminars. 

\begin{theorem}\label{main} Let $K_0$ and $K_1$ be  $2$-foams without boundary embedded in $4$-space with diagrams $D_i$ for $i=0,1.$ 
$K_0$ and $K_1$ are isotopic  if and only if there is a sequence $D_{j/n}$   of diagrams for $j=0,\ldots,n$  such that $D_{j/n}$ differs from $D_{(j-1)/n}$ by one of the Reidemeister/Roseman-type moves that are listed below.
\end{theorem}

\newpage 
\begin{center}
\includegraphics[width=6.5in]{allmoves}
\end{center}

The movie parametrizations do not always agree with the projections as drawn, but they are topologically equivalent. In the right-most column not all necessary crossings are indicated, but these are easy to guess given the remaining information. The proof of Theorem~\ref{main} will be presented in Section~\ref{prmain}. First (Section~\ref{dim3}), we discuss the  Reidemeister moves for knotted trivalent graphs. Then (Section~\ref{camel}) we develop a short digression on the interactions among critical points, vertices and crossings. These critical interactions provide the atomic pieces used to construct foams. A second digression (Section~\ref{Turaev}) develops an idea first presented by Turaev~\cite{TuInven} that indicates why only certain moves are needed as long as all possible type-II moves are present.  Section~\ref{cp4v} is a short section on critical points of foams and their edges. The idea of the proof of Theorem~\ref{main} is to analyze the possible critical points and intersections among the self-intersections and the singular strata. Section~\ref{wrapup} examines two additional moves that change the topology of a foam, but for which a regular neighborhood is preserved. Section~\ref{future} points to unfinished endeavors.

\newpage

\section{Critical Points and Crossings --- $1$-Dimensional Case}
\label{dim3}

\begin{wrapfigure}{r}{3.2in}\vspace{-.2in}\includegraphics[width=3in]{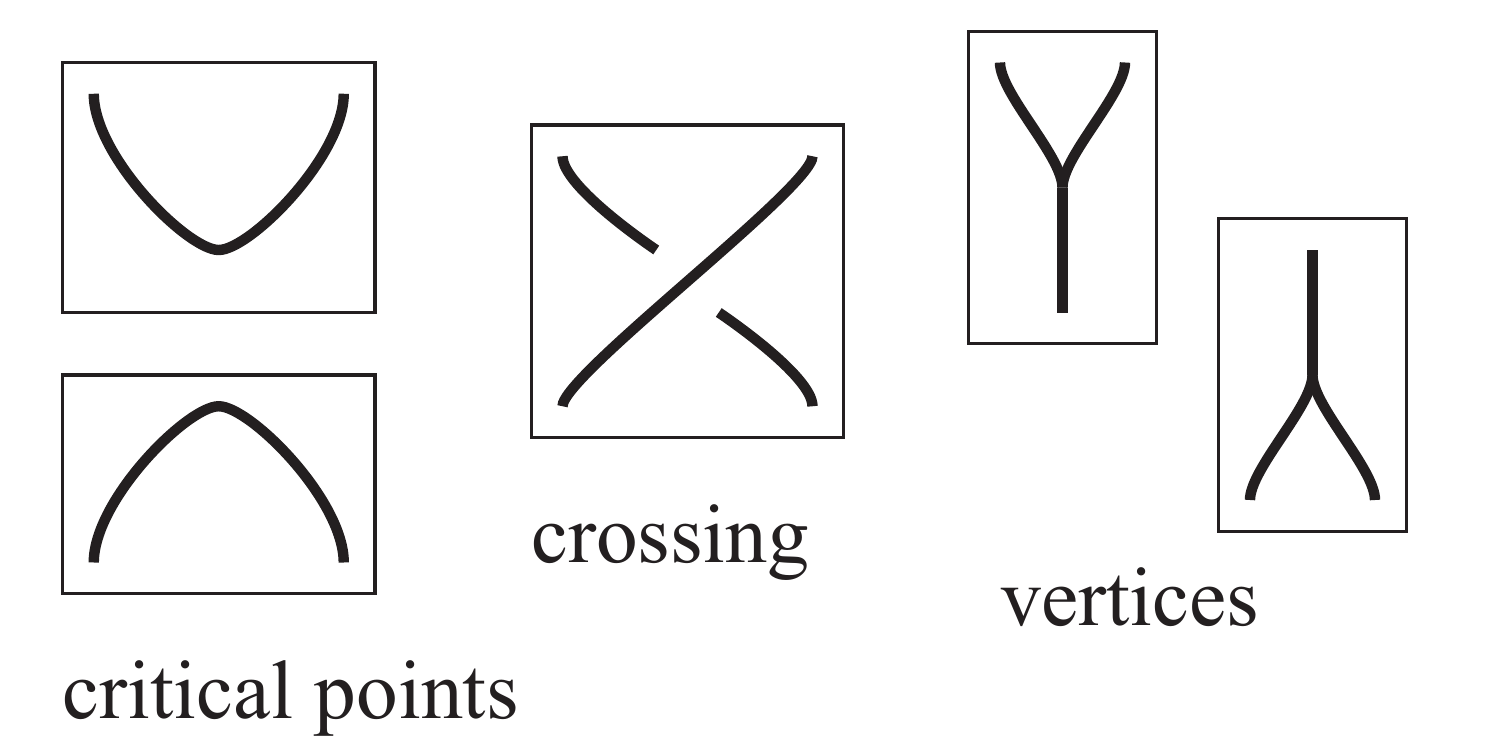}\end{wrapfigure}
Consider a trivalent graph that is embedded in $3$-dimensional space. A generic projection onto the plane will have isolated transverse double points and points of no higher multiplicity. The points to which trivalent vertices project are not double points. If a height function is chosen in the plane, the graph may be assumed to have non-degenerate critical points that are either maxima or minima.

\begin{wrapfigure}{l}{1.5in}\includegraphics[width=1.3in]{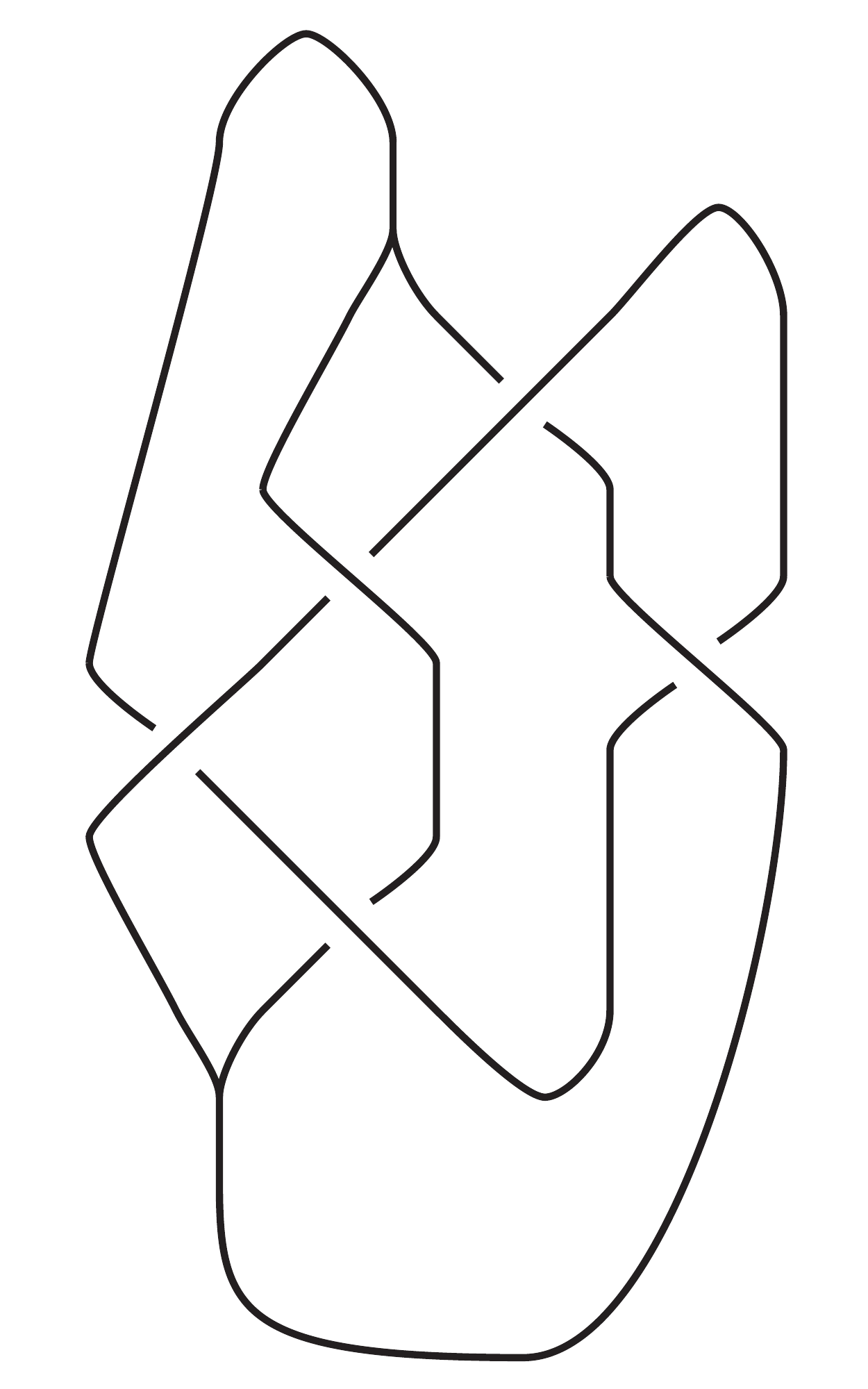}\end{wrapfigure}
 Furthermore, one can rearrange the graph in space so that the crossing points, vertices, and critical points all lie at distinct levels as indicated to the left.

The vertices of {\sf Y}, the crossing points, and the critical points are $0$-dimensional. To quantify the moves to trivalent graphs, we examine the transverse intersections and critical points of the corresponding sets in ${\R^2} \times [0,1]$ as an isotopy occurs. For example, a Reidemeister type-II move is a critical point of the $1$-dimensional crossing set that is engendered as the projection of spacial graph moves in the plane.
 A Reidemeister type-III move is the transverse intersection between the trace of a crossing and the $2$-dimensional sheet consisting of an arc of the diagram times the isotopy parameter. Similarly, there are two scenarios in which a vertex passes through a transverse arc: in one the vertex passes below the arc while in the other the vertex passes under the arc. A Reidemeister type-I move corresponds to a critical point on the double decker set during the isotopy. The twisting of a trivalent vertex is analogous to a Reidemeister type-I move. 
 
 \begin{wrapfigure}{r}{2.7in}\vspace{-.2in}\includegraphics[width=2.5in]{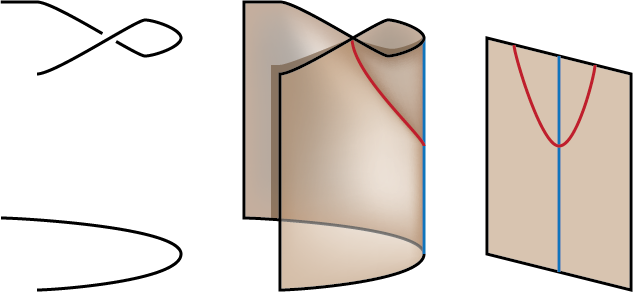}\vspace{-.2in}\end{wrapfigure} 
 In the figures that flank this paragraph and those that follow, the Reidemeister-type moves for trivalent spacial graphs and the corresponding critical events in the surfaces that represent the isotopies are illustrated. The written discussion will now elaborate upon the nature of critical points and transverse intersections.  
 
 First consider a type-I move (indicated as an RI-move). A single crossing is involved. The crossing occurs along an arc of the knot diagram. As the loop surrounded by the repeated crossing shrinks, the double-decker points on the arc converge to a simple critical point. This move, then, corresponds to a generic critical point of the double-decker set. 
  
  \begin{wrapfigure}{l}{2.7in}\vspace{-.2in}\includegraphics[width=2.5in]{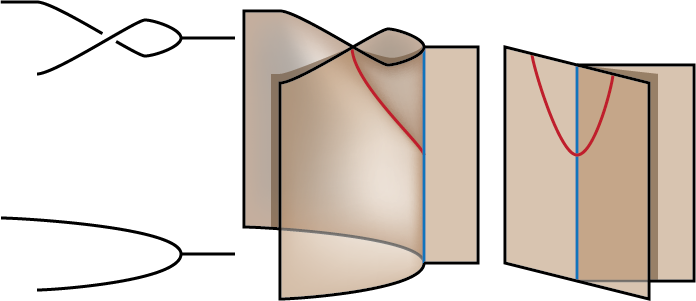}\end{wrapfigure} The {\it twisted vertex move} (Tw), that is depicted to the left, also corresponds to a critical point of the $1$-dimensional double decker set   of the isotopy. We remark here that in the case of both the twisted vertex move and the type-I move, only one type of crossing is illustrated. Obviously, the move also holds with the opposite crossing.   
 
  \begin{wrapfigure}[7]{r}{2.7in}\vspace{-.5in}\includegraphics[width=2.5in]{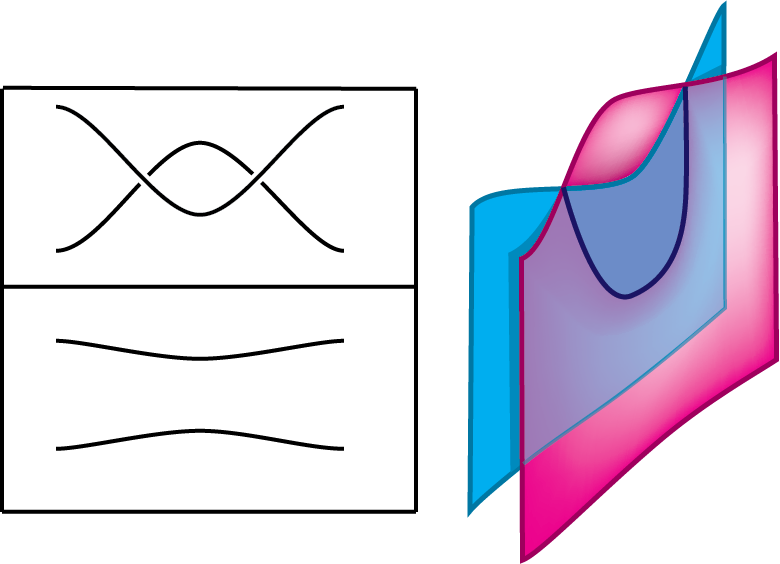}\end{wrapfigure} The {\it type-II move} (RII) is indicated to the right. 
It involves a pair of distinct double points that converge during the isotopy to a critical point on the double point set. Note that on the double-decker set there are a pair of critical points --- one for each sheet involved in the crossing --- and the critical levels coincide since they are equivariant with respect to the involution on the double-decker set.   
 
The double points and the vertices of a trivalent graph are the $0$-dimensional singularities. The critical points of the double point set, then correspond to the type-I and type-II Reidemeister moves, and the twisted vertex. Critical points for the set of vertices do not strictly correspond to moves for trivalent spacial graphs since they affect the topology of the graph. They will, however, represent aspects of the corresponding foam, and they will be discussed in Section~\ref{cp4v}.

 The {\it type-III move} (RIII) and the moves in which a vertex passes over (IY) or under  (YI) a transverse arc are all manifestations of a $1$-dimensional set passing transversely through a $2$-dimensional set in the $3$-dimensional space of the isotopy direction times the plane of projection. The moves are depicted as broken surface diagrams below the current paragraph. 

\includegraphics[width=6in]{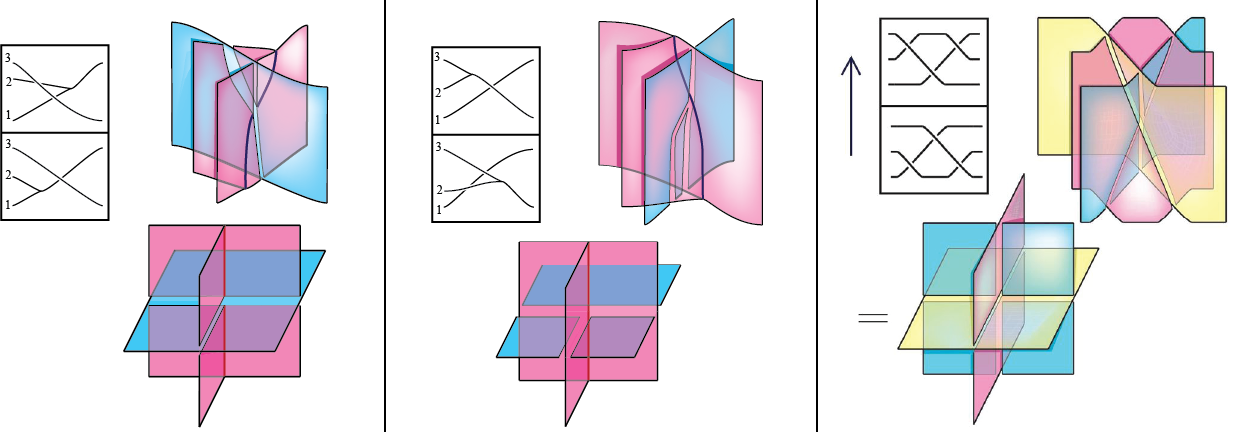}

  \begin{wrapfigure}{l}{2.7in}\vspace{-.3in}\includegraphics[width=2.5in]{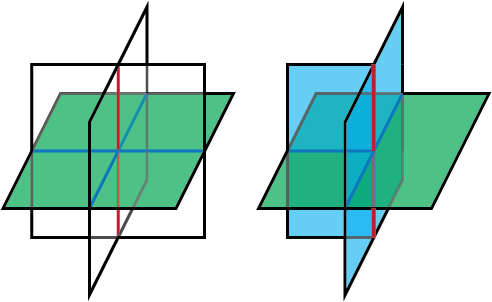}\end{wrapfigure} To the left of the current paragraph, the projections of the type-III move and either move in which a vertex passes through a transverse arc are depicted. The $1$-dimensional set is a vertical arc that is indicated in red, and the $2$-dimensional set is horizontal and depicted in green. In the case of the triple point, any one of the crossings  between top/middle, top/bottom, or middle/bottom could serve as a red arc while either the bottom arc, middle arc, or the top arc (respectively) traces out the green sheet. In the case of {\sf Y}$\times [0,1]$ intersecting the transverse sheet, the arc formed from the vertex intersects the transverse sheet which is either entirely above or entirely below the {\sf Y} $\times [0,1]$.  

Next we complete the proof of the following classical result. See for example \cite{AB}, \cite{Reid}.

\begin{theorem}\label{cl} For $i=0,1$, let $k_i: G \hookrightarrow \R^3$ be spacial embeddings of a trivalent graph $G$ represented by diagrams $D_i$.  Then $k_0(G)$ and $k_1(G)$ are isotopic embeddings if and only if there is a sequence of diagrams $D_{j/n}$ for $j=0, \ldots, n$ such that $D_{j/n}$ differs from $D_{(j-1)/n}$ by a planar isotopy or an application of one of the moves indicated in the figure  below.  \end{theorem}
\begin{center}
\includegraphics[width=2in]{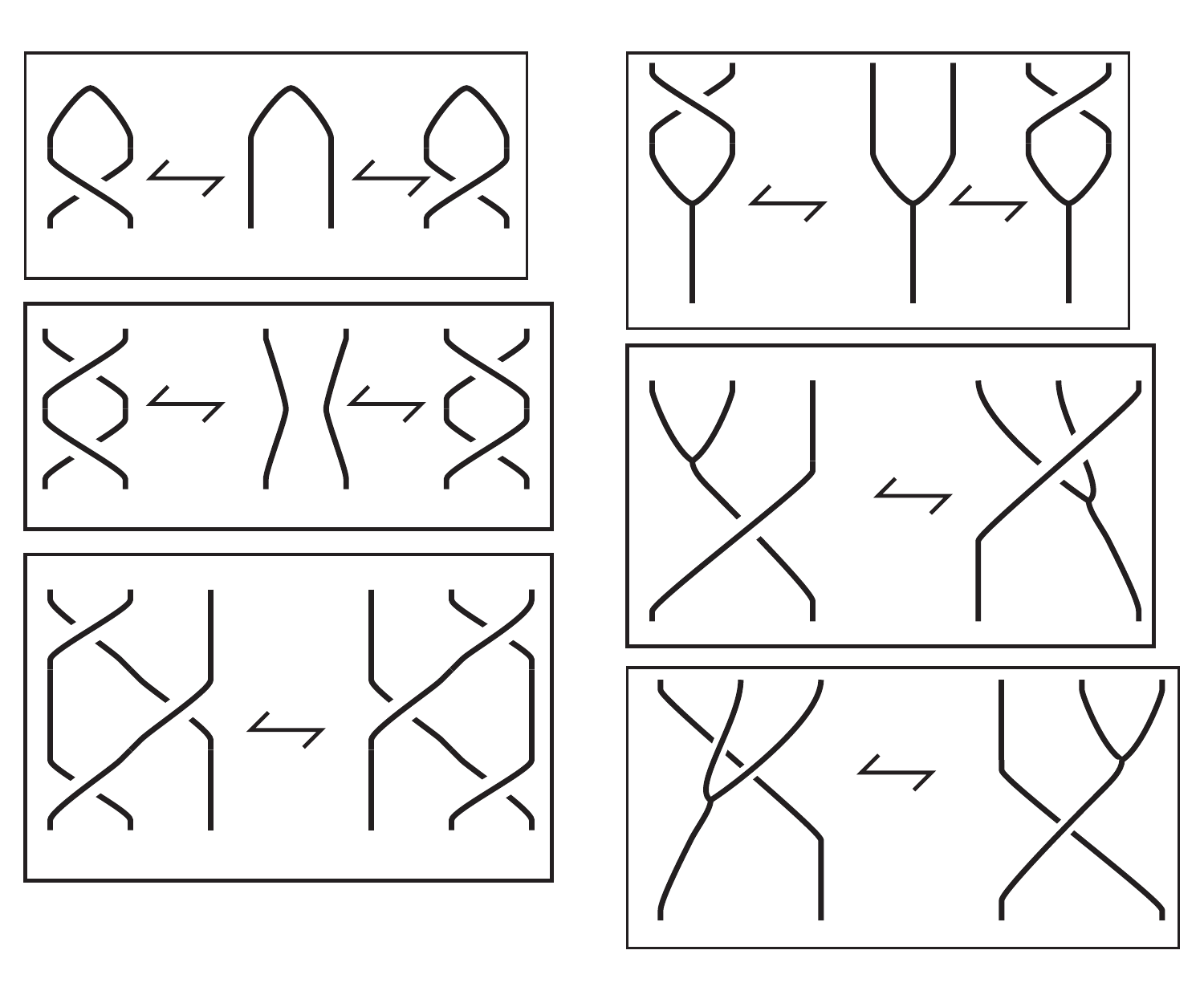}
\end{center}
{\sc Proof.} If the diagrams $D_0$ and $D_1$ differ by a finite sequence of these moves, then they are (clearly) isotopic. 

Suppose that the two diagrams are isotopic. 
The $0$-dimensional sets for the diagram $D_0$ are the crossing points and vertices. During an isotopy, these trace out $1$-dimensional sets. Using the isotopy direction as a Morse function for this $1$-dimensional set, we quantify the critical points. As mentioned above, the critical points of the vertex set alter the topology of the underlying graph. So we only consider the critical points of  crossing set. These correspond to type-I and type-II Reidemeister moves. The transverse intersections between an arc of crossing points (or edge of ${\mbox{\sf Y}}  \times[0,1]$) and the $2$-dimensional sheets formed from arcs of the diagram times the isotopy parameter correspond to the type-III moves or the vertex of a {\sf Y} crossing over or under a transverse sheet (the YI or IY moves). 

Thus a given isotopy may be modified so that the critical points and transverse intersections occur at differing times. Moreover, by compactness, we may assume that there are only finitely many such critical points or intersections. This completes the proof. $\Box$

We remind the reader that the move considered in the second row, second column of Theorem~\ref{cl} is called {\it the YI-move}. The move in the third row, second column is called {\it the IY-move}. The move in the first row of the second column is called {\it twisted vertex move}.

\section{Critical Points, Crossings, and Vertices}\label{camel}

\begin{wrapfigure}[14]{r}{2.5in}\vspace{-.5in}\includegraphics[width=2.4in]{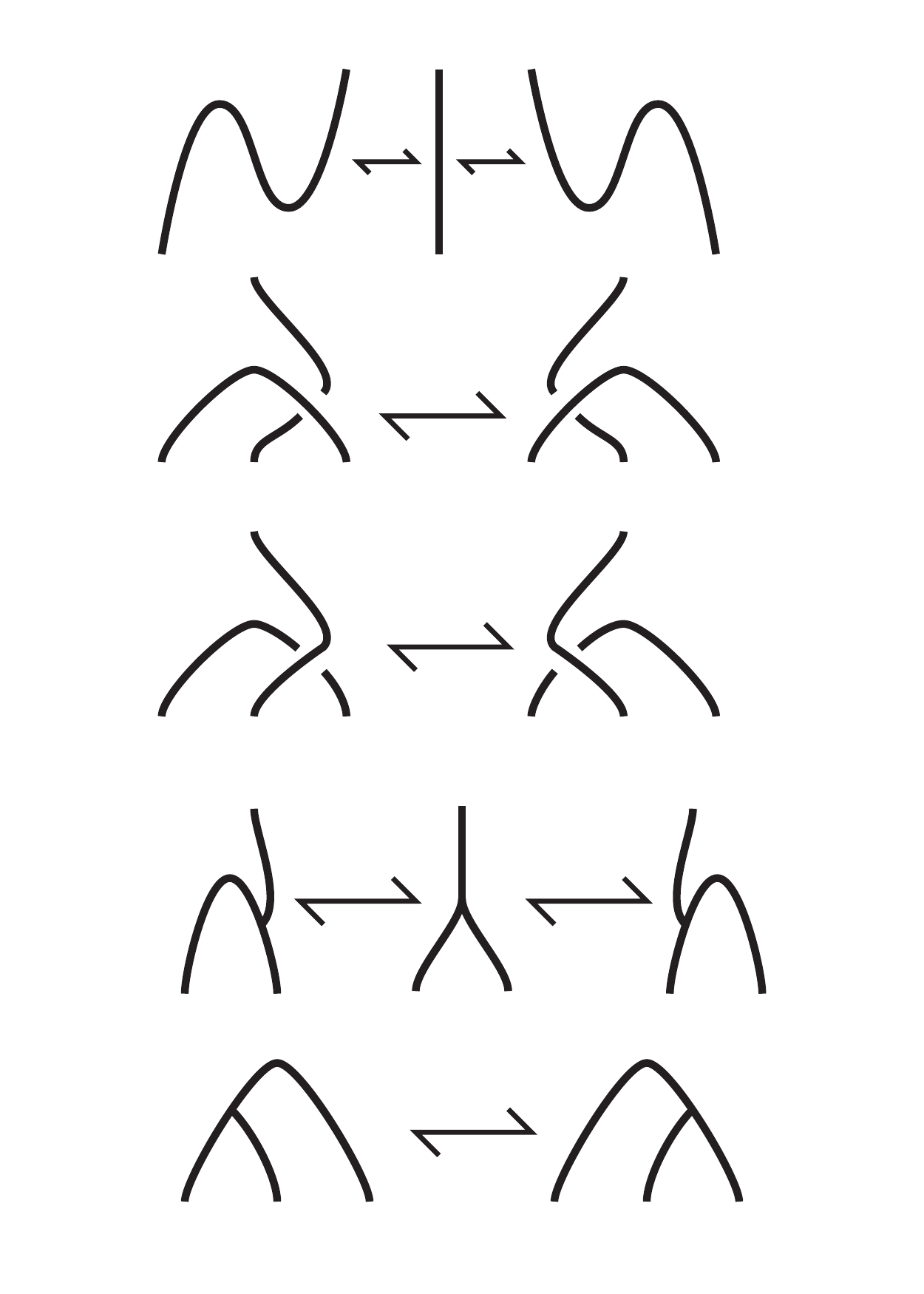}\end{wrapfigure}
In many circumstances, it is necessary to impose a height function upon the plane into which a knotted trivalent graph is projected. For example, in developing categorical/algebraic interpretations of the graph, a height function is chosen, the trivalent vertex represents a multiplication or comultiplication operator, and the generic moves that reflect the changes in height are axioms in the algebraic setting. In particular, multiplication and comultiplication can be defined in terms of each other in the Frobenius algebra setting. This coincidence between algebraic structures occurs as a result of generic perturbations of the height functions.   The moves indicated to the right together with the exchange of distant critical points are sufficient to re-adjust the height function of any embedded trivalent graph.

\begin{wrapfigure}[15]{l}{2.5in}\includegraphics[width=2.4in]{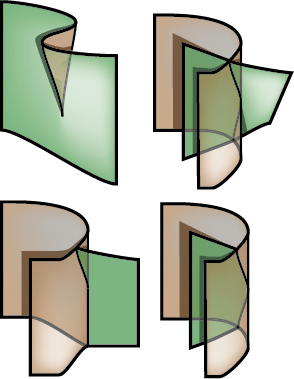}\end{wrapfigure} Each of these moves can be interpreted as a local picture of the projection of a knotted $2$-foam in $4$-space into the plane. In this case, a foam has generic folds that either end in cusps or end along the edge of a foam. Thus a fold ends in a cusp,  a double point curve can pass over a fold,  an edge passes over a fold, or a fold is created at the trivalent edge. In this paper, we will not study the interactions of such moves, nor will we attempt to incorporate them into a movie-move theorem such as that given in \cite{CRS} or \cite{R2B} which contains a heuristic analysis for immersed surfaces.

Instead, the interpretations are presented here since they facilitate the illustration of knotted foams. Specifically, we can define a diagram of a knotted foam as a projection with crossing information indicated along the double curves. When an author desires to illustrate a foam, he or she can make a detailed movie in which each still in the movie has a height function, and slow the action down until successive stills differ by simple moves. The moves are birth, death, zippers (which are all illustrated below), type-I, type-II, type-III, IY, YI,  associators, and the commutation of distant critical points). Then the author can use the local projections to draw the pieces of a knotted foam. The drawing is, of course, in the plane of the paper, and keeping track of the folds helps indicate many of the details of the foam.

It is an interesting aspect of the geometry that the corresponding moves with minima in place of maxima follow from the moves given here. In particular the zig-zig moves (at the top of the previous illustration) and the commutation of distant critical points allow the ``upside-down" versions to hold. While this seems to be well-known among experts, and indeed it is a standard exercise, a proof is provided below.  
\begin{center}\includegraphics[width=4in]{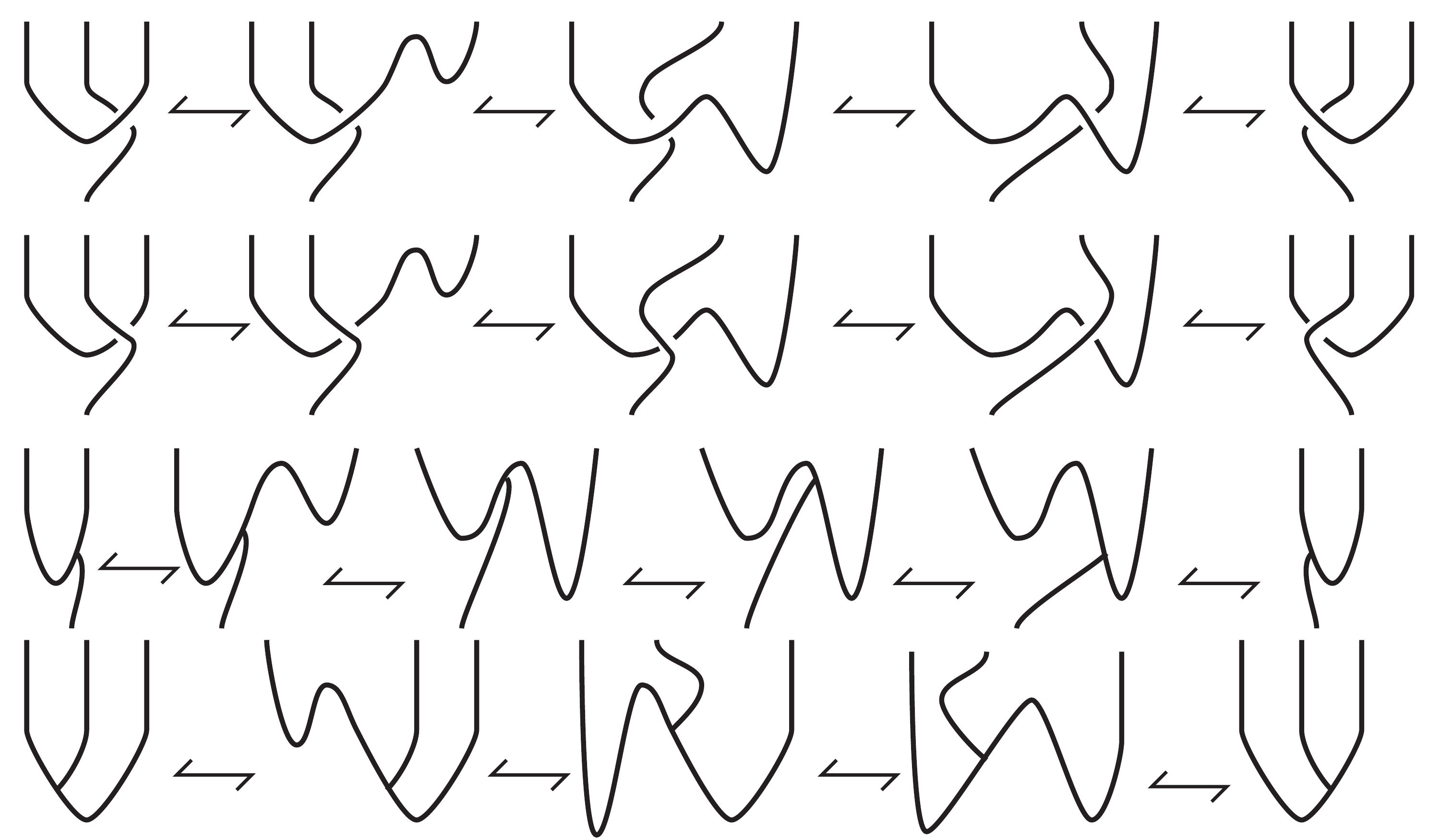}\end{center}

\section{Turaev's trick}\label{Turaev}

In this section, we continue to work in the un-oriented category. For a type-III move there are six possible initial crossings for which the move is valid. Each of these follows from the type-III move that is given and the two possible type-II moves given. Similarly, the YI and IY moves are stated with specific crossing information given. Using Turaev's idea, we indicate that the moves with alternate crossing information given follow from the given moves and the type-II moves given. It is not difficult to generalize these proofs to the situations in which arcs are oriented. In the case of the Reidemeister type-III moves, the situation can be thought of very algebraically. Starting from the braid relation $aba=bab$, we prove: $\overline{a}ba$ $=$ $\overline{a}bab \overline{b}$ $=$ $\overline{a}aba \overline{b}$ $=$ $ba \overline{b}.$ Similarly, $ab\overline{a}$ $=$ $\overline{b} bab\overline{a}$ 
 $=$ $\overline{b} aba \overline{a}$ $=$ $\overline{b} ab$. The remaining three identities follow similarly. 
\begin{center}\includegraphics[width=6.5in]{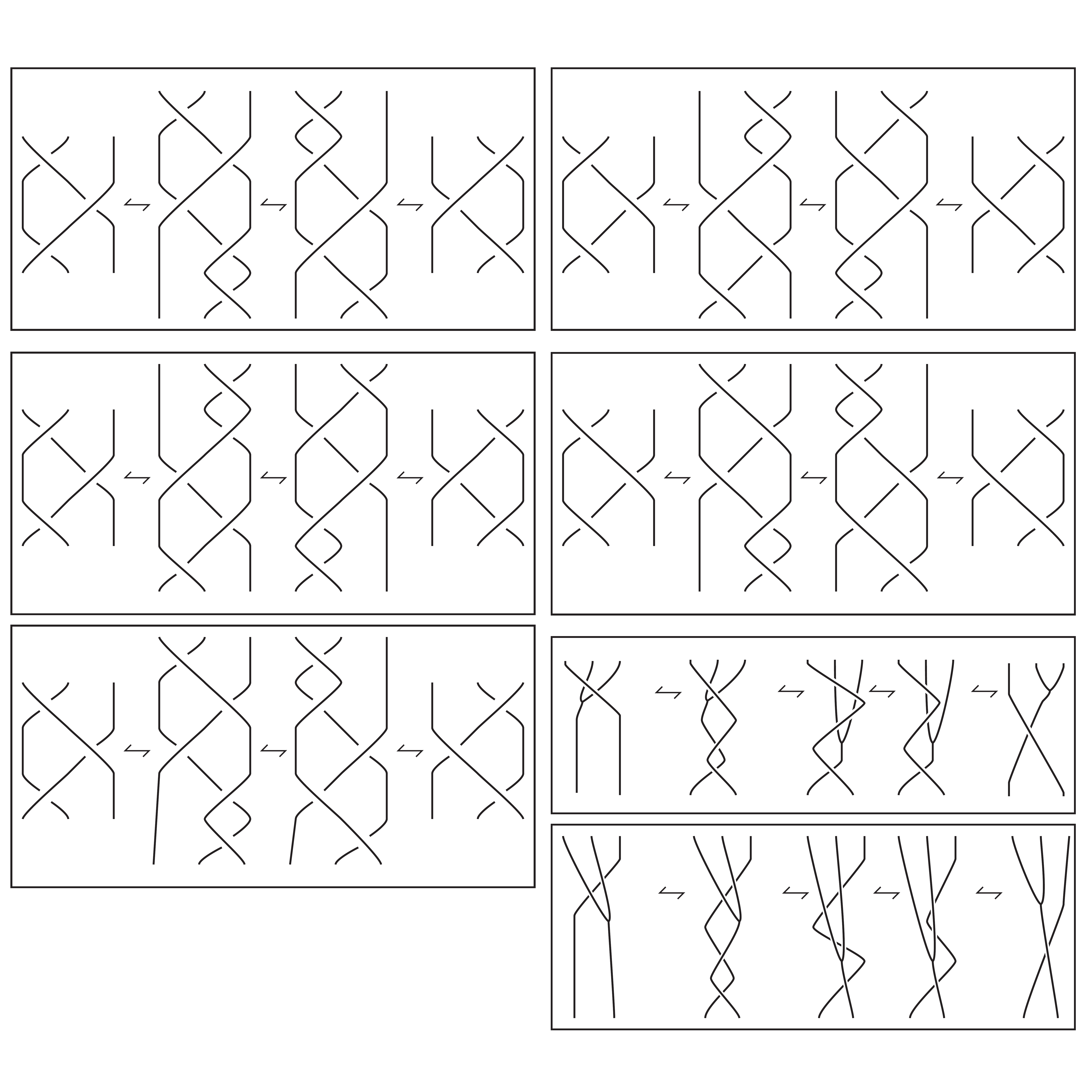}\end{center}

\newpage

\section{Critical points for foams}
\label{cp4v}

\begin{wrapfigure}[6]{r}{3in}\vspace{-.25in}\includegraphics[width=2.75in]{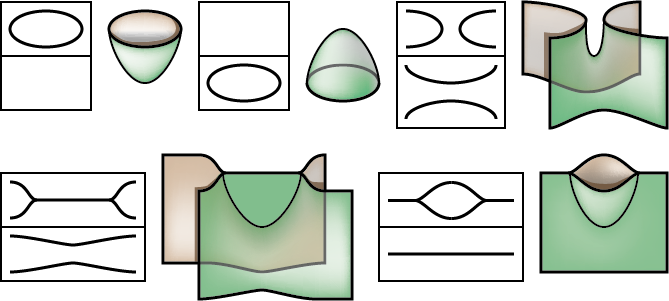}\end{wrapfigure}
Let $k: F \hookrightarrow \R^4$ denote an embedding of a closed foam $F$ into $4$-space. Consider a generic projection $p: \R^4 \rightarrow \R^3$. Then after a small perturbation, if necessary, a height function on the image $p(k(F))$ can be chosen so that the critical points of the foam, its edge set, and the double point set are all at distinct levels. Moreover, the generic triple points, the intersections of transverse sheets with the edge set, the vertices are all at distinct levels, and these levels are distinct from the critical points of the previous sentence. The critical points are births or deaths of simple closed curves or saddle points as indicated to the right. The critical points for the edge set are also illustrated here; these are called {\it zipper moves}.  

The remaining local pictures are the time-elapsed versions of the type-I,II, and III Reidemeister moves, the YI and IY moves, the associator move (which is the neighborhood of a vertex in a foam), and a twisted vertex move (Tw). As mentioned above, in drawing (projecting to the plane), the edges and double curves can pass over the folds which correspond to the traces of the critical points in the stills. 

The analysis of this section allows us to describe a knotted foam by means of a movie when necessary.

\section{Critical and intersection behaviors for isotopies of foams}
\label{prmain}

We turn now to proving the main result, Theorem~\ref{main}. The sketch of the proof goes as follows. Branch points, twisted vertices, triple points, and the intersections between a transverse sheet and the edges of the foam are $0$-dimensional and hence isolated. The non-degenerate critical points 
in the isotopy directions correspond to each of the moves RI, Tw, RIII, YI, and IY being invertible. For the RI, Tw, YI, and IY moves, there are two types of invertibilty: elliptic and hyperbolic, that depend on the structure nearby. The critical points of the double point set correspond to the Roseman bubble and saddle moves for the double point set. These correspond to the elliptic and hyperbolic confluence of the double point set or the two types of invertibility of the the RII-move. 

The next set of moves occur when, in the isotopy direction, a transverse sheet intersects any one of the $0$-dimensional sets. That is, an embedded sheet becomes an embedded $3$-dimensional solid in the $4$-dimensional space-time of the isotopy direction. Meanwhile, branch points, twisted vertices,  YI and IY intersections,  triple points, and foam vertices evolve in the time direction and yield embedded arcs in space-time. The transverse intersection between an arc and a $3$-dimensional solid in $4$-space is an isolated point. These transverse intersections account for five of the remaining moves. 

Finally, we can consider the double points or the edges of foams  evolving in space-time to become $2$-dimensional sets. The transverse intersection between two surfaces in $4$-space consists of isolated points. These transverse intersections can also be used to reinterpret some of the previous moves; more importantly, they include a move below that is called the YY-move. The rest of the proof consists of recognizing that all of the possible critical points and transverse intersections have been identified.

\begin{wrapfigure}[21]{r}{4in}\vspace{-.25in}\includegraphics[width=3.78in]{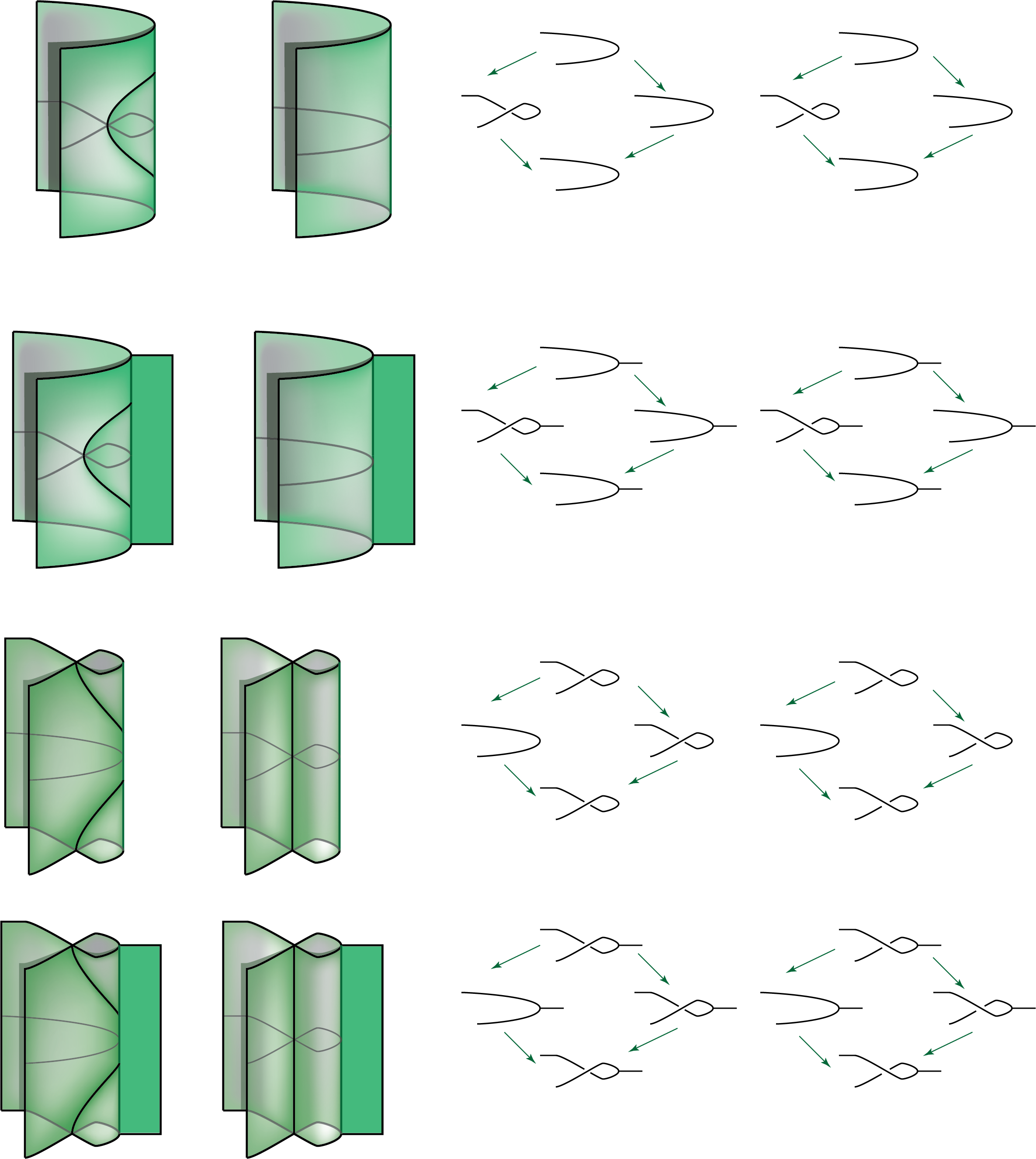} \end{wrapfigure}
Four moves to foams are illustrated to the right of this paragraph that indicate the critical points of the branch points and twist vertices. The movie versions indicate the two possible liftings into $4$-dimensional space. The branch points and twisted vertices are $0$-dimensional and therefore are isolated. In the isotopy direction, these points form $1$-dimensional sets  whose non-degenerate critical points correspond to the moves indicated. At the top of the illustration, an {\it elliptic confluence of branch points} is indicated. Analogously, and immediately below is an {\it elliptic confluence of twisted vertices}. The next two moves are the {\it hyperbolic confluence of branch points or of twisted vertices}, respectively. In the movie parametrizations, we consider the source and target graphs to appear as the top or bottom boundary of either foam. The path around the left indicates the foam on the left and that on the right indicates the foam on the right. It is easy to imagine a singular graph that encapsulates the singularity that occurs as the branch points or twisted vertices converge.

\begin{wrapfigure}[9]{l}{2.5in}\includegraphics[width=2.35in]{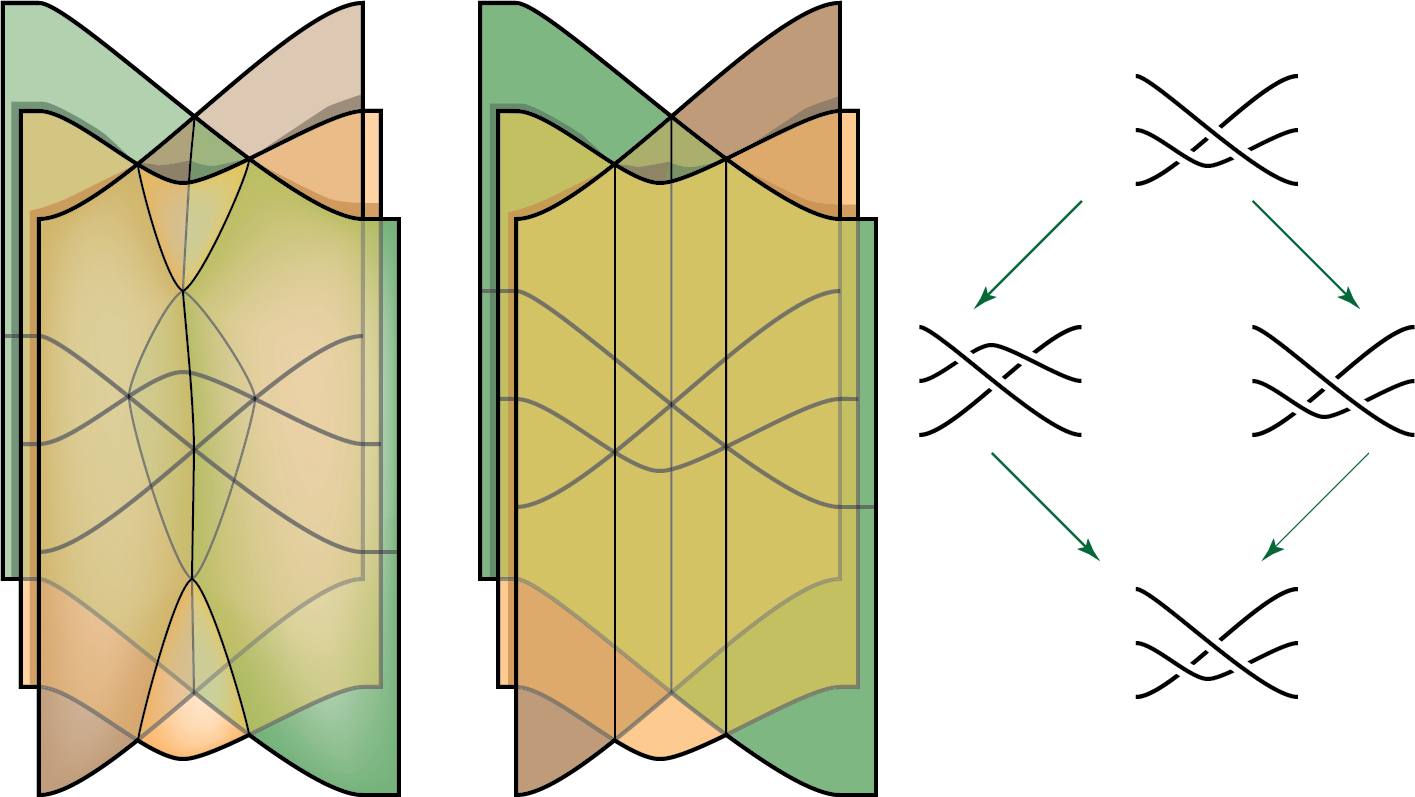}  \end{wrapfigure}
In a similar fashion, the triple points of the projection of a knotted foam are $0$-dimensional and isolated. In the isotopy direction, a non-degenerate critical point corresponds to the annihilation/creation of a pair triple points. Only one of six possible lifts is indicated in the movie version here. However, all six possible liftings are  needed in order to facilitate the higher dimensional analogue of Tureav's trick. The move illustrated here is called {\it type-III type-III-inverse move}.

\begin{wrapfigure}{r}{3in}\vspace{-.2in}
\includegraphics[width= 2.75in]{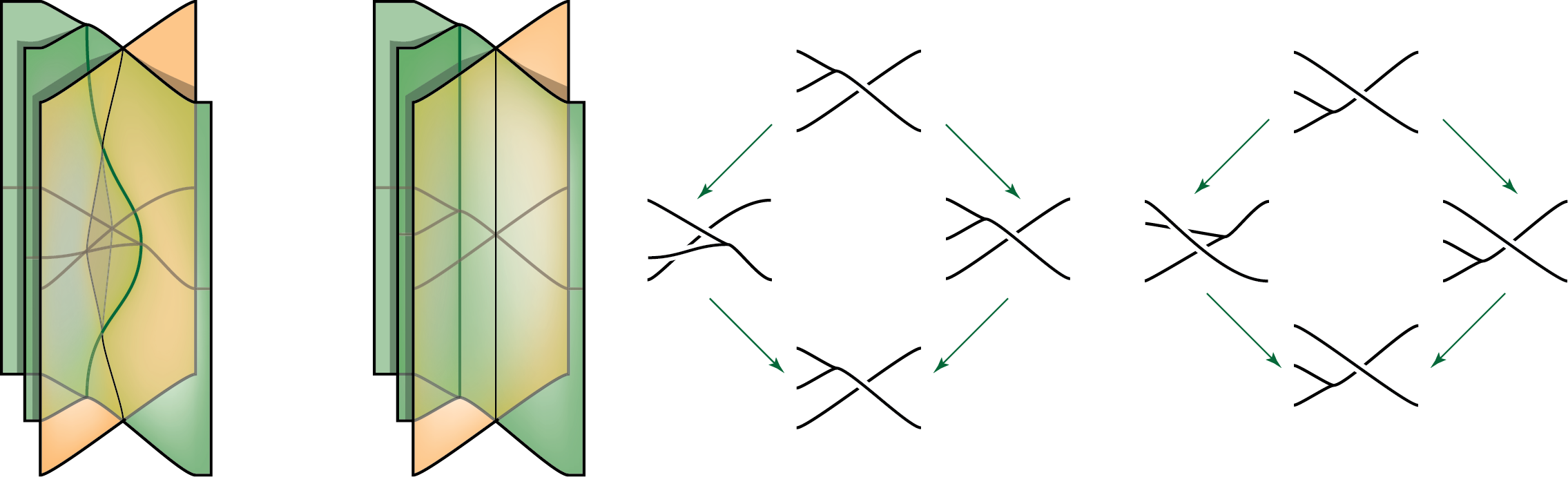} 
\end{wrapfigure} This move is a {\it  YI-bubble move} or a {\it  IY-bubble move}. The crossing points between an edge and a transverse sheet are $0$-dimensional. A critical point in the isotopy direction creates or annihilates such a pair. There are two possible sets of crossings for such a move. These correspond to the YI or the IY moves as illustrated to the right of the projected surfaces.

\begin{wrapfigure}[7]{l}{3in}
\includegraphics[width= 2.75in]{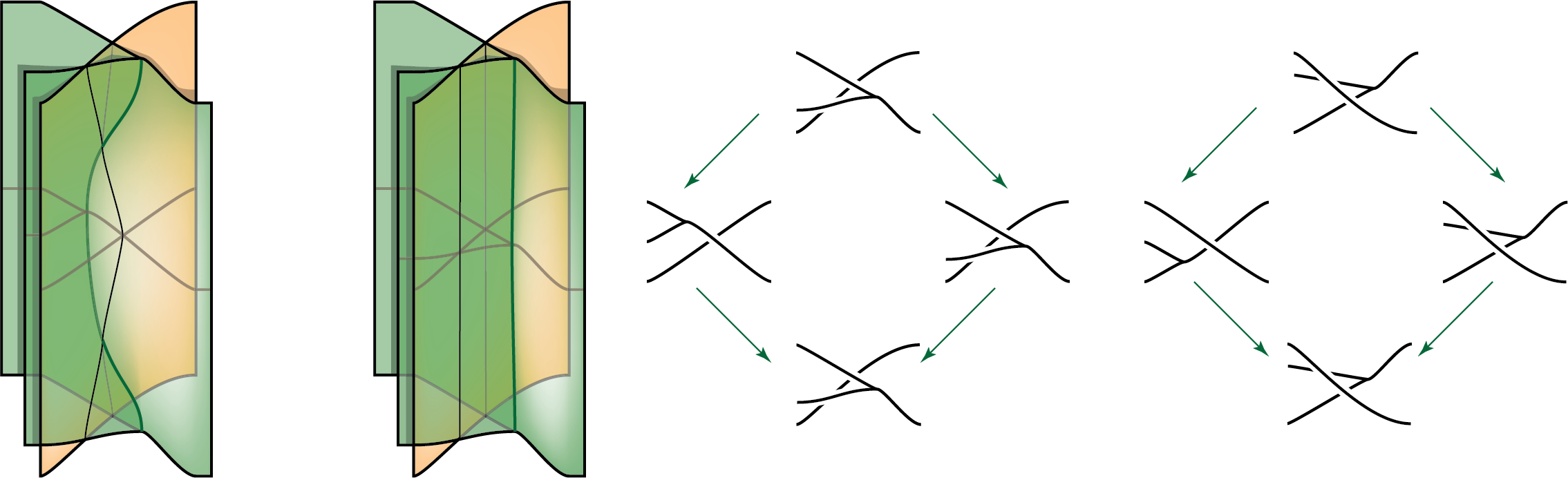} 
\end{wrapfigure} The move to the left is  a {\it YI-saddle move} or an {\it IY-saddle move}.  The analysis of the crossing points and the edge of the foam follows exactly as in the preceding paragraph. On the left-hand-side of the move, there are a pair of $0$-dimensional multiple points that are formed as the transverse intersection between an edge of the foam and an embedded sheet. A critical point in the isotopy direction annihilates or creates such a pair of crossing points. At the critical point, the edge of the foam is tangent to the embedded sheet. There are two lifts into $4$-space of this move and these are represented by the movie moves to the right of the illustrations of the foams.

\begin{wrapfigure}{r}{3.5in}
\includegraphics[width= 3.25in]{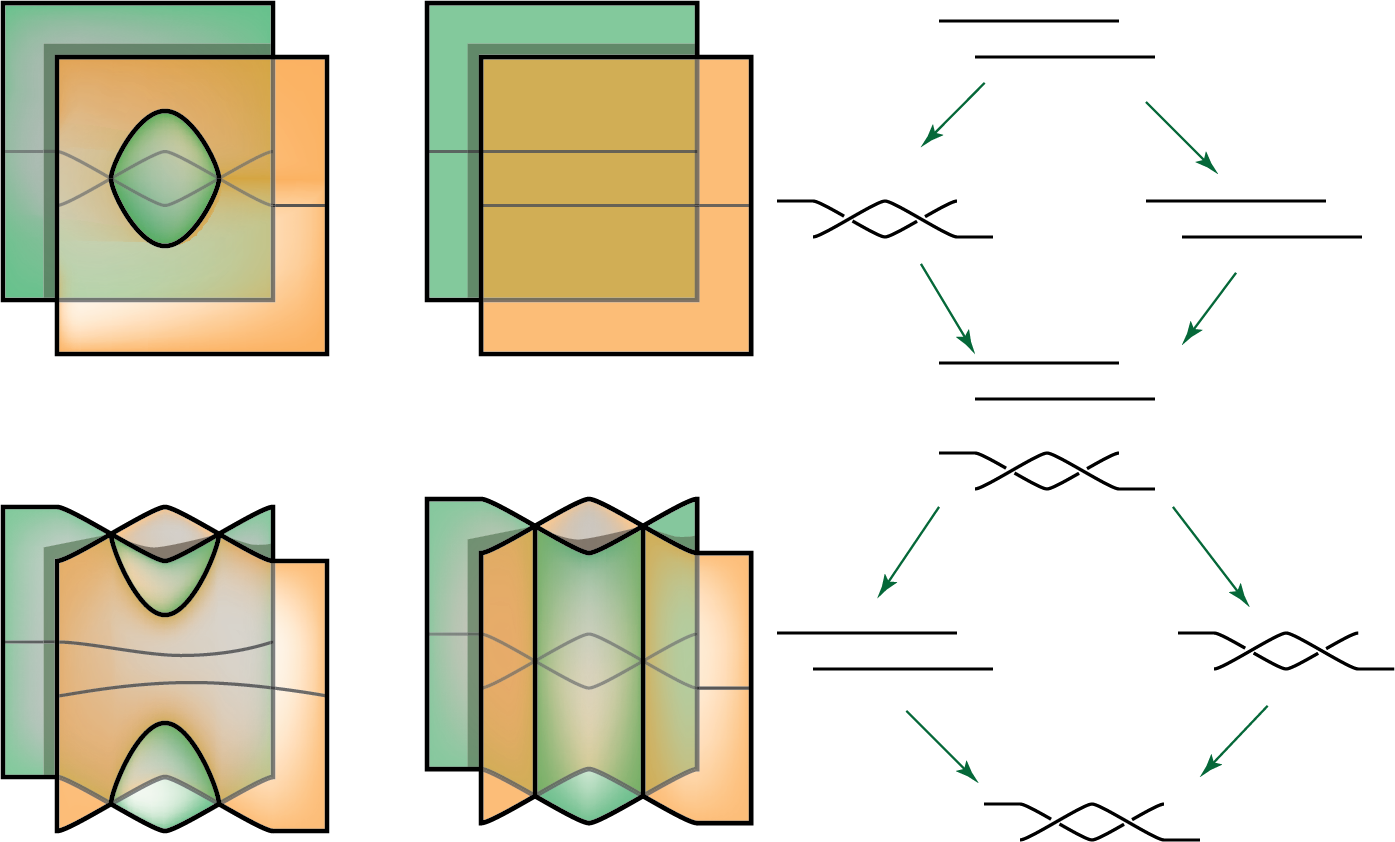} 
\end{wrapfigure}  The moves illustrated to the the right are the {\it type-II bubble} or {\it type-II saddle moves}. These correspond to optimal or saddle critical points of the $1$-dimensional double point set of the foam. Both are found in Roseman's list of moves to knotted surfaces since neither involves an edge of a foam. In the isotopy direction the double point set is $2$-dimensional, and these moves represent surface critical points. In the movie parametrizations only one of the two possible crossings are illustrated.

\begin{wrapfigure}{l}{3.5in}
\includegraphics[width= 3.25in]{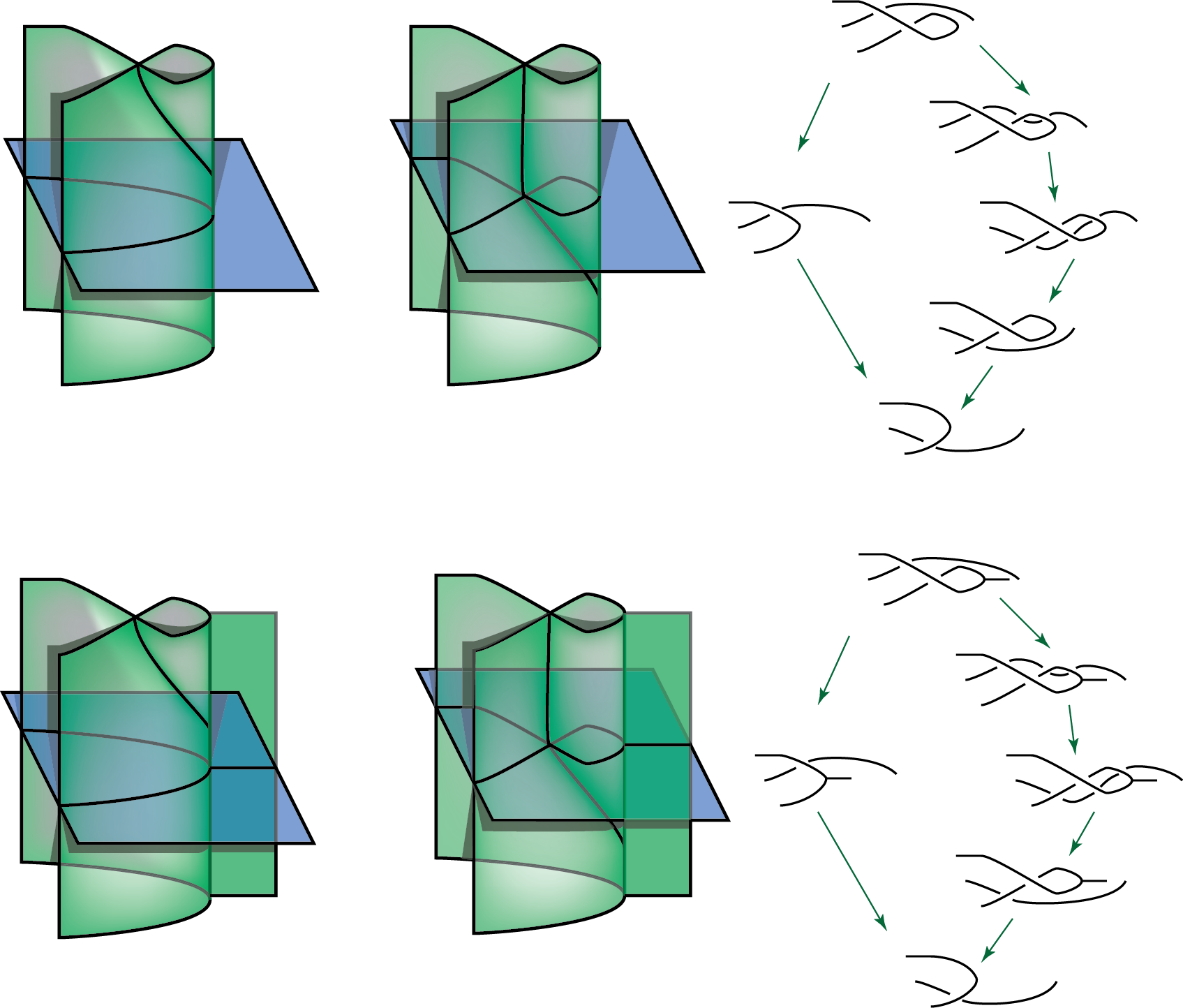} 
\end{wrapfigure} The two illustrations to the left indicate {\it pushing a branch point (or a twisted vertex) through a transverse sheet.} The transverse intersection between a $3$-dimensional solid and a $1$-dimensional arc in $4$-space is an isolated point. The branch point or twist vertex evolves in space-time to be an arc and the nearby embedded sheet evolves to be $3$-dimensional. As a result, on the right-hand-side of either move there is a triple point. The movie parametrizations do not indicate both possible liftings. Furthermore neither coincides with a strict interpretation as horizontal cross sections of the figures: it is more convenient to draw the transverse sheet horizontally. Nevertheless, there is a sequence of slices that gives the movie versions indicated.

The three movie moves that are indicated below the current paragraph are called from left to right, {\it the IYI-move}, {\it the IIY-move}, and the {\it YII-move}. When these three moves are combined with all possible type-III type-III-inverse moves, then any foam of the form ${\mbox{\sf Y}}  \times[0,1]$ can pass over, under, or through a pair of transverse sheets. To make that statement more precise, let us consider the YII-move that is depicted upon the right.  The ${\mbox{\sf Y}}  \times [0,1]$ is completely below the two transverse sheets. Another possibility is that it lies above the two transverse sheets. Or it could be above the second sheet and below the third sheet. In any of these cases the movie move can be performed. Similar situations hold for the IYI and the IIY-moves. However, in the presence of a sufficient class of type-III type-III inverse moves, these alternative crossings follow from those that are indicated here.
\begin{center}
\includegraphics[width= 6in]{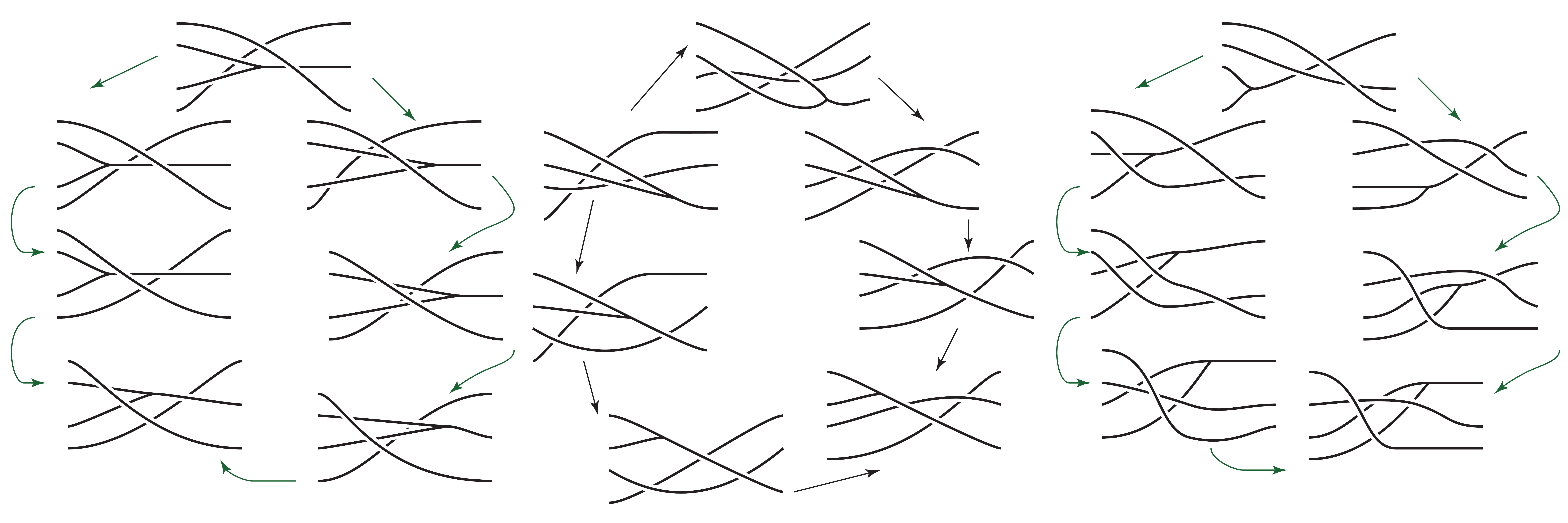} 
\end{center}

\begin{wrapfigure}{r}{3.5in}
\includegraphics[width= 3.25in]{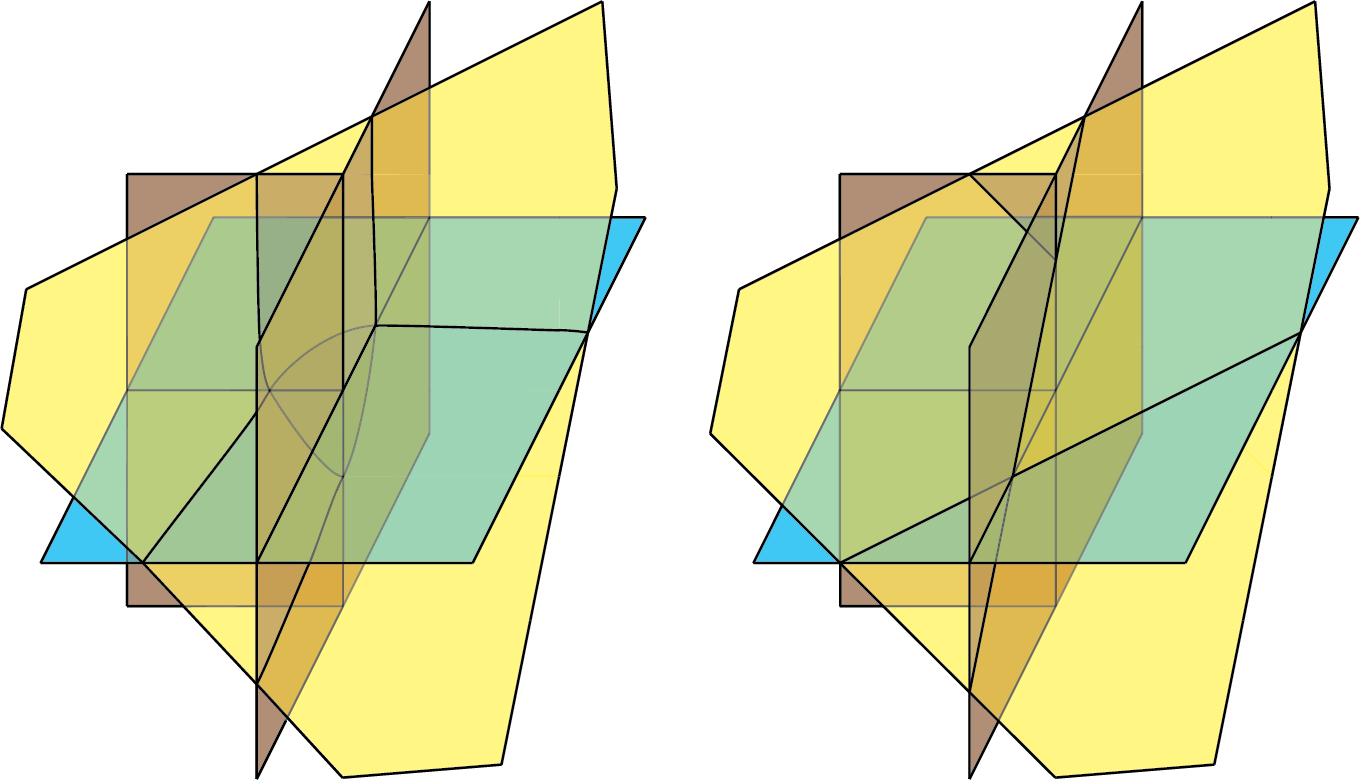} 
\end{wrapfigure} 
Each of the YII, IYI, and IIY moves projects in space-time to the same configuration. On the right,  a foam of the form ${\mbox{\sf Y}}  \times [0,1]$ intersecting a pair of intersecting embedded disks. We can consider the intersection of an edge of ${\mbox{\sf Y}}  \times [0,1]$ with one of these sheets as an isolated vertex that can pass through the remaining sheet. (In space-time, this is a $1$-dimensional sheet --- crossing $\times$ interval --- intersecting a $3$-dimensional sheet). Alternatively, we may consider the edge of ${\mbox{\sf Y}}   \times [0,1]$  moving through the double line of the remaining two sheets. In this case the transverse intersection in space-time is between two $2$-dimensional surfaces. Either analysis of the transverse intersections gives rise to the move.

\rule{6.5in}{0in}
\begin{wrapfigure}{l}{4in}
\includegraphics[width= 3.85in]{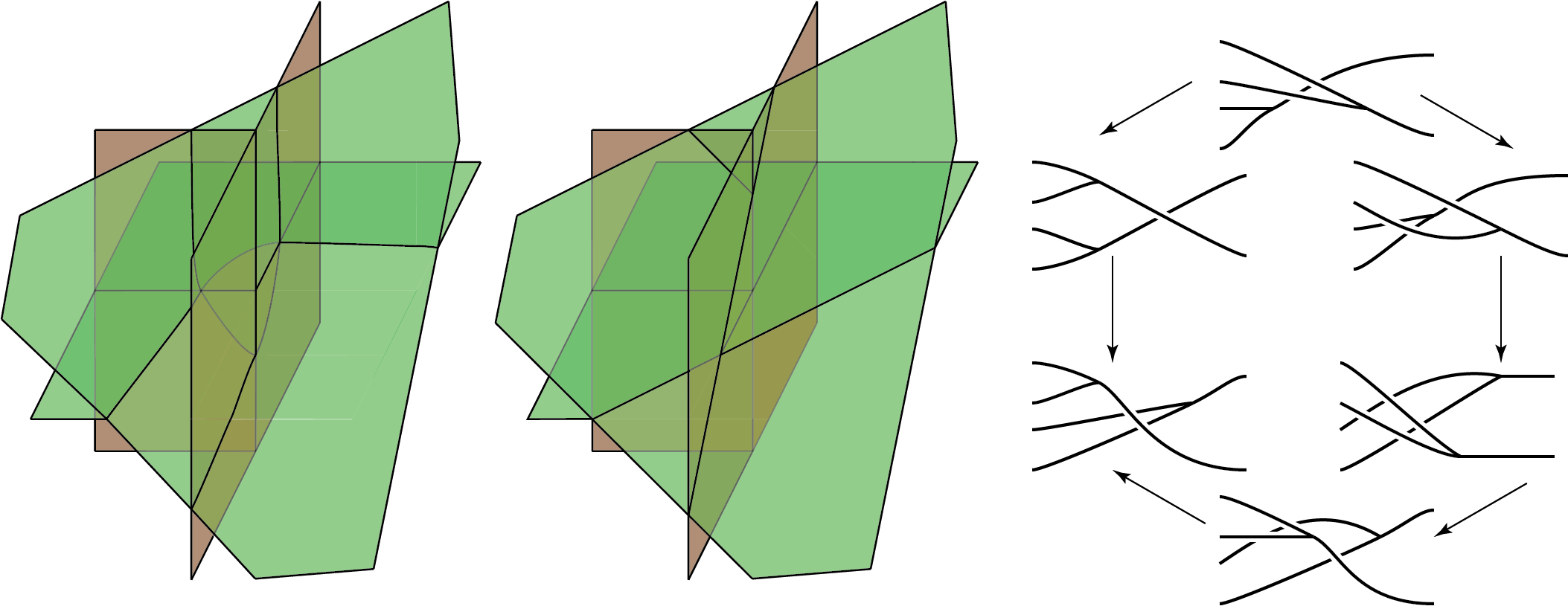}
\end{wrapfigure}  {\it The YY-move} is illustrated here. Two edges, each coming from a foam of the form ${\mbox{\sf Y}}  \times [0,1]$,  pass through each other in space-time.  In the isotopy direction, an embedded edge evolves to be an embedded disk in space-time. The transverse intersection of two $2$-disks in $4$-space is an isolated point. The movie illustration indicates crossing information. The alternative crossing information can be obtained by reflection. There is not a precise match between the movie parametrization and the surfaces as projected here, but it is not difficult to rearrange the surfaces by isotopy in $3$-space to make the figures match.

\rule{6.5in}{0in}
\begin{wrapfigure}{r}{4in}
\includegraphics[width= 3.75in]{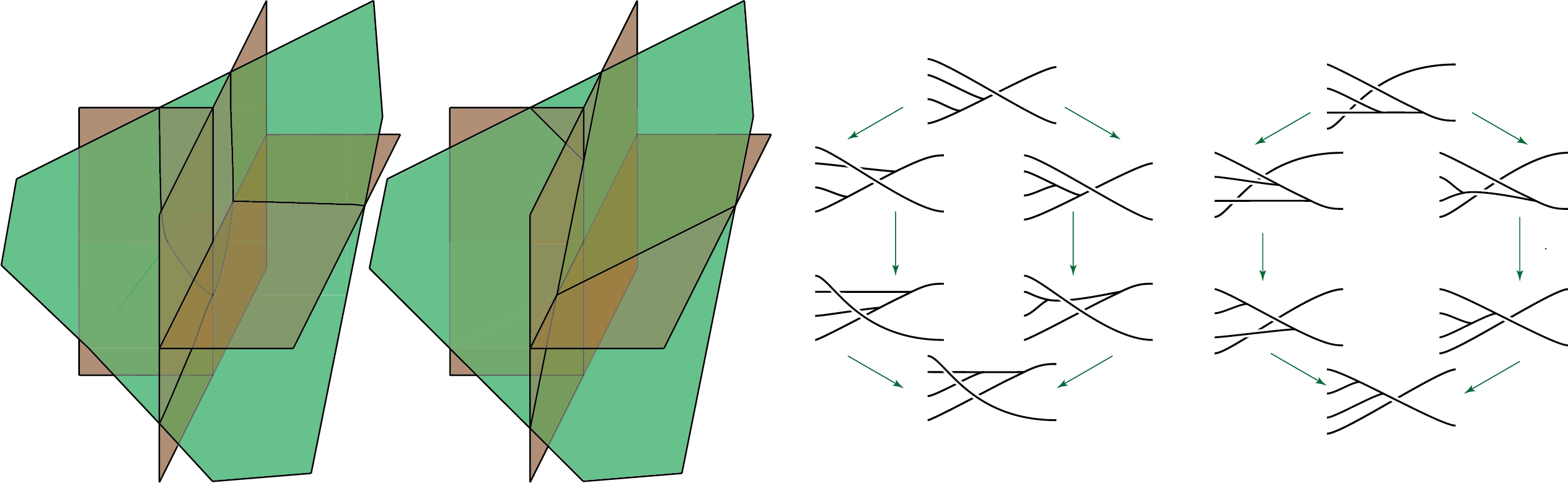} 
\end{wrapfigure} 
{\it The YYI-move} and the {\it IYY-move} are depicted here. 
A vertex at the juncture of the foam $Y^2$ passes through a transverse sheet. In space-time the vertex evolves into a $1$-dimensional set, and as before, the transverse embedded sheet  becomes $3$-dimensional. The two movie moves that are illustrated include crossing information. Either the transverse sheet is completely above or completely below the foam $Y^2$.

\rule{6.5in}{0in}
\begin{wrapfigure}{l}{4in}
\includegraphics[width= 3.75in]{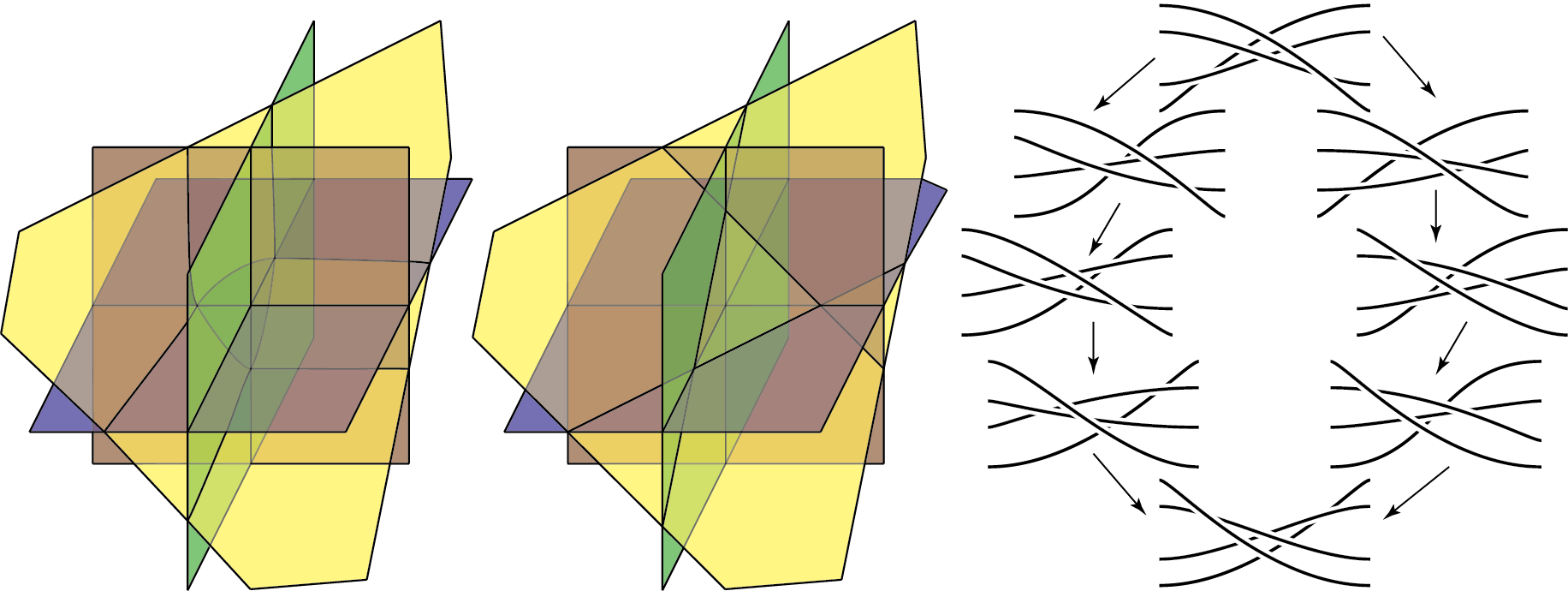} 
\end{wrapfigure} 
{\it The tetrahedral move, or quadruple point move} is indicated here with one possible  collection of crossing information. The moves that involve different crossing information follow from a higher dimensional analogue of Turaev's trick and the various type-III type-III-inverse moves (see \cite{CJKLS} for an implementation). One of the triple points (chosen among four possible) passes through a transverse sheet. The triple point evolves in space-time to form a $1$-dimensional arc. The transverse sheet evolves as a $3$-dimensional solid. The transverse intersection in space-time is an isolated point. This point is also the intersection of two double point arcs. There are three ways to parse the double point set into a pair of disjoint arcs. These are counted as pairs of two element subsets of $\{1,2,3,4\}$; specifically, $\{ \{1\cap 2\}, \{3\cap 4\}\},$ $\{ \{1\cap 3\}, \{2\cap 4\}\}$, or $\{ \{1\cap 4\}, \{2\cap 3\}\}$ where the numbers indicate  labels on each of the surfaces that intersect on either side of the move.

In the eleven paragraphs that precede the current paragraph, each of the potential Roseman-type moves has been described as a critical point of a multiple point stratum or the transverse intersection between  intersection strata or embedded sheets. We observe that these descriptions are exhaustive. The $0$-dimensional sets are the branch points, twist vertices, triple points, or intersections between edges of the foam and a transverse sheet. The $1$-dimensional sets are the double point arcs. Clearly, the $2$-dimensional sets are the embedded sheets in the foam. Every possible critical point or transverse intersection has been accounted. Thus these are the codimension $1$ singularities.


Now let us complete the proof of Theorem~\ref{main}. First, if two diagrams differ by any one of the moves depicted, then they are isotopic. 
Next suppose that two diagrams represent isotopic embeddings. We have a map ${\mathcal K}: F \times [0,1] \rightarrow \R^4$ such that the restrictions ${\mathcal K}|_{F \times \{i\}} = K_i$ to the ends represent the given knottings of the foam $F$. Moreover, for each $t\in [0,1]$, the foam $F\times \{t\}$ is embedded. The isotopy ${\mathcal K}$ can be adjusted slightly, if necessary, so that the composition $p_2 \circ {\mathcal K}$ is has generic singularities on the singular sets of $3$-dimensional foam that is the product $F \times [0,1]$. By compactness, there will be finitely many singularities. The isotopy can be perturbed further, if necessary, so that each singular point lies at a different time coordinate. The types of generic singularities that need to be quantified are critical points for the $1$ and $2$-dimensional strata of the isotopy and transverse intersections between strata. 

The catalogue of $0$-dimensional singular points on a given knotted foam are as follows:
\begin{enumerate}
\item[(0.1)] branch points or twist vertices that result from an RI or Tw move to a knotted trivalent graph,
\item[(0.2)] triple points or the intersection between an edge of ${\sf Y} \times [0,1]$ and an embedded sheet, or
\item[(0.3)] vertices of the foam.
\end{enumerate}
The $1$-dimensional singular points of a knotted foam are as follows:
\begin{enumerate}
\item[(1.1)] the double point arcs that are caused by the transverse intersection between a pair of sheets of the foam, or
\item[(1.2)]  the edge set of a foam.
\end{enumerate}
Any face of the foam $F$ is $2$-dimensional and therefore non-singular.

Under the assumption that the singularities of the isotopy are generic and isolated, we examine each singular situation. 

First we consider the critical points of $1$-dimensional sets. The singular points described in items (0.1), (0.2), and (0.3) evolve in space-time to $1$-dimensional sets. 
During the isotopy, the topology of the foam remains constant. Therefore,  there are no critical points for the vertex set listed as item (0.3). (However, in the following section, we discuss such topological changes).
The critical points for branch points and twist vertices (0.2) are the elliptic and hyperbolic confluences of each. These are accounted for in the first illustration of this section.  The critical points for the RIII, YI, or IY  moves correspond to the type-III type-III inverse moves, the YI or IY bubble  and saddle moves. 

In examining the critical points of the $2$-dimensional sets in the isotopy, we first observe (in a manner similar to that above) that the edge set of the foam $F$ remains unchanged during the isotopy. Thus we only have to consider the critical points of the double point set. These correspond to the type-II saddle and bubble moves.

The remaining singularities to be considered are those that are caused by the transverse intersections between two strata. In $4$-space, we have intersections between a $3$-dimensional stratum and a $1$-dimensional stratum (3-1), or between a pair of $2$-dimensional strata (2-2). To track the (3-1) intersections, we consider arcs formed in space-time by the vertices listed in items (0.1) through (0.3) and examine their transverse intersections with an embedded sheet. The resulting intersections are as follows:
\begin{enumerate}
\item[(3.1.B)] a branch point passing through a transverse sheet,
\item[(3.1.Tw)] a twisted vertex passing through a transverse sheet,
\item[(3.1.V)] a vertex passing through a transverse sheet as in the YYI move,
\item[(3.1.YI)]  the intersection between an edge and a transverse sheet passing through a third sheet (moves YII,IYI, and IIY), 
\item[(3.1.III)] the intersection between a triple point and a fourth transverse sheet (quadruple point move),
\item[(2.2.YY)] the intersection between two edges of the foam (move YY),
\item[(2.2.YD)] the intersection between an edge of the foam and a double curve (moves YII, IYI, and IIY), and
\item[(2.2.DD)] the intersection between two double point arcs (the quadruple point move).
\end{enumerate}
We observe that there are coincidences between (3.1.III) and (2.2.DD), as well as coincidences between (3.1.YI) and (2.2.YD). 

Finally, we point out that the assumptions on general position are warranted. First, the isotopy can be fixed within a neighborhood of its $1$-dimensional strata (the trace of the vertices, branch points, twist points, triple points, and edge intersections with transverse sheets) so that the critical points of these sets are non-degenerate. Next, without affecting the $1$-dimensional strata,  the double point set can be adjusted to have non-degenerate critical points. Any adjustment to the double point set can be extended to a neighborhood of the double point set within the foam. Intersections between a pair of double point arcs, a pair of edges or a double point arc can also be made to be generic and transverse as can intersections between vertices and sheets. There perturbations are then extended over the remainder of the $F \times [0,1]$. This completes the proof.

\section{Moves that change the topology of the underlying foam}
\label{wrapup}

It is important to remark that not only are embedded trivalent graphs studied for their own sake, but a given knotted handlebody in $3$-space deformation retracts to an embedded trivalent graph. The graph, however, is not unique. Two graphs that ``carry" such a knotted handlebody differ by the so-called IH-move. Up to equivalence, the IH-move is  given via the movie parametrization of the basic foam $Y^2$. The theory of knotted handlebodies embedded in $3$-space is equivalent to the theory of knotted trivalent graphs modulo the IH-move. 

\begin{wrapfigure}{r}{2.75in}
\includegraphics[width= 2.5in]{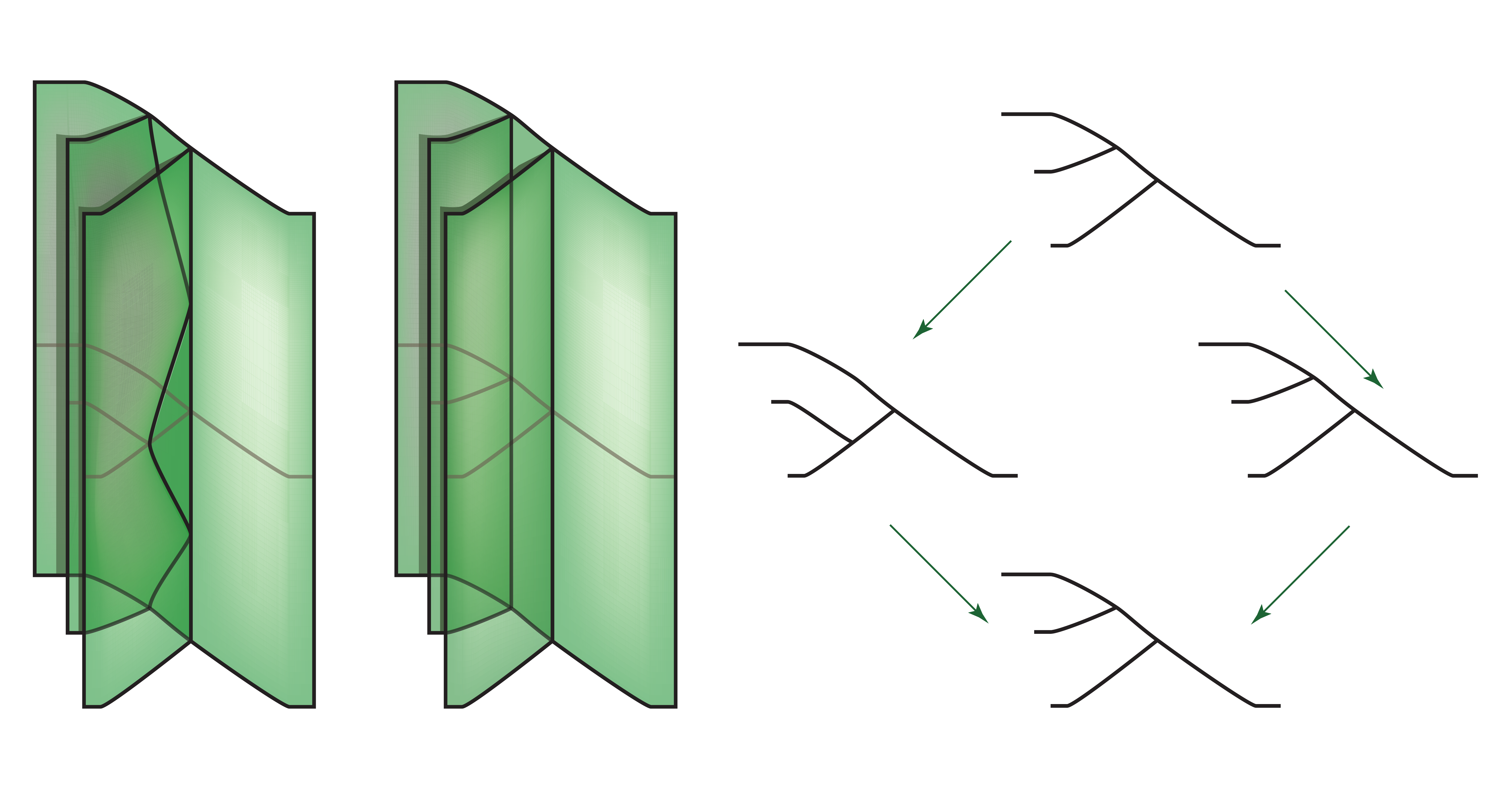} 
\end{wrapfigure} A similar situation holds in $4$-space. We can include among the Roseman moves two additional moves that are indicated below.  The first of these is the invertibility of the IH-move. In the theory of special spines for $3$-manifolds it is sometimes called the {\it lune move} or the {\it orthogonality condition}. The reason for the latter name comes from the Tureav-Viro \cite{TV} invariants, the neighborhood of a vertex of a foam is colored by representations of $U_q(sl_2)$, and the move corresponds to the orthogonality condition for the $6j$-symbol. See also~\cite{CFS}.
Observe that if a foam is embedded in $3$-space, then regular neighborhood of the foam is invariant under this move. So similarly, a regular neighborhood in $4$-space of such a foam is also invariant since it can be obtained from the neighborhood in $3$-space by the cartesian product with an open interval.

\begin{center}
\includegraphics[width= 3.75in]{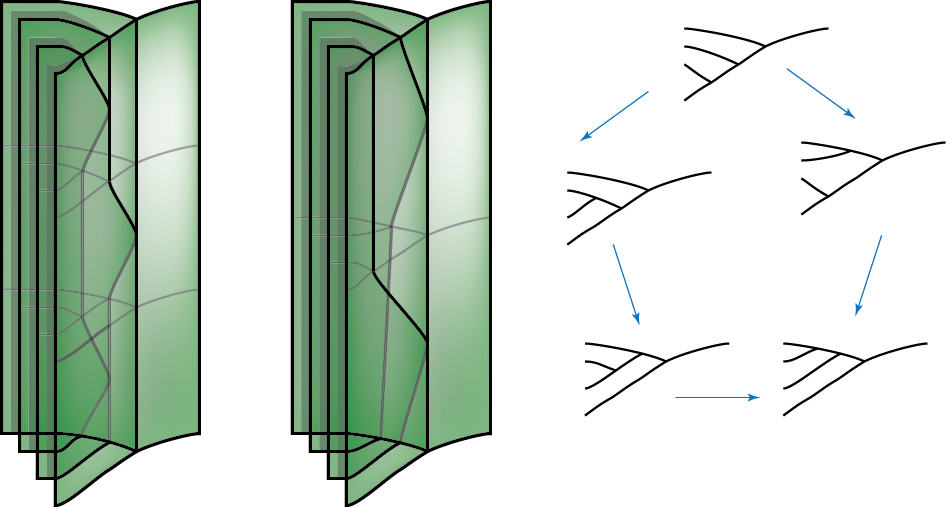} 
\end{center}

The remaining move is the {\it (3,2)-move} or {\it Elliott-Beidenharn} move. Again when considered as a move to special spines, following \cite{Matveev,Piergallini}, the move preserves the topology of the underlying $3$-manifold. It is easy to see that a regular neighborhood of the foam is preserved under the move in $3$-space and hence in $4$-space. The resulting theory is the theory of embedded $4$-manifolds in $4$-space. 

The lune move corresponds to a critical point of the vertex set in deformation of a foam, and the $(3,2)$-move corresponds  
to a vertex in a $3$-dimensional foam. 
Thus it is reasonable to study isotopy classes of knotted foams modulo these additional two relations. 

\section{Future work}\label{future}

This paper is a technical piece that is necessary for a serious study of knotted foams and their $3$ (and higher)-dimensional generalizations. In work with Atsushi Ishii and Masahico Saito, we will establish a cohomology theory for certain algebraic systems that is sufficient to define non-trivial invariants of knotted $2$-foams in $4$-space. 

The inclusion of the penultimate section is also meant to indicate the initial stages in the study of $3$-dimensional foams. In particular, one can construct movies of $3$-dimensional foams embedded in $5$-space by including the Roseman/Reidemeister moves of Theorem~\ref{main}, the orthogonality and Eilliott-Beidenharn moves, critical points of surfaces and critical points of the edge sets. Thus the critical points of surfaces correspond to $0$, $1$, $2$, and $3$-handles that are attached to the solid sheets of $3$-foams. Moreover the edge set  of a $3$-foam is a $2$-foam. Important moves to $3$-foams are easy to establish. A full Roseman-type theorem is unknown to the author at this time, but determining one should be routine and geometrically tedious. 

Finally, it is worth mentioning that there is an underlying categorical motivation here that is related to the tangle hypothesis of Baez and Dolan \cite{BD}. Here we are considering the interaction between a braiding and a Frobenius structure as well as the identities among relations of these. The precise location of knotted foams in the Baez-Dolan table is an interesting taxonomic problem.


\begin{flushleft}
J. Scott Carter \\ 
Department of Mathematics \\ 
University of South Alabama \\ 
Mobile, AL 36688 
USA \\
E-mail address: {\tt carter@southalabama.edu} 
\end{flushleft}

\end{document}